\documentclass[reqno,12pt]{amsart}
\usepackage{amsfonts}
\usepackage{amsmath,amssymb,amsthm}

\newtheorem{thm}{Theorem}[section]

\newtheorem{lem}[thm]{Lemma}

\theoremstyle{definition}

\theoremstyle{remark}

\begin{document}
\title[MERIDIAN SURFACES WITH POINTWISE 1-TYPE GAUSS MAP]{MERIDIAN SURFACES
IN $\mathbb{E}^{4}$ WITH POINTWISE 1-TYPE GAUSS MAP}
\author{Kadri Arslan, Bet\"{u}l Bulca and Velichka Milousheva}
\address{Uluda\u{g} University, Department of Mathematics, 16059 Bursa,
Turkey}
\email{arslan@uludag.edu.tr}
\address{Uluda\u{g} University, Department of Mathematics, 16059 Bursa,
Turkey}
\email{bbulca@uludag.edu.tr}
\address{Bulgarian Academy of Sciences, Institute of Mathematics and
Informatics, Acad. G. Bonchev Str. bl. 8, 1113, Sofia, Bulgaria; "L.
Karavelov" Civil Engineering Higher School, 175 Suhodolska Str., 1373 Sofia,
Bulgaria}
\email{vmil@math.bas.bg}

\begin{abstract}
In the present article we study a special class of surfaces in the
four-dimensional Euclidean space, which are one-parameter systems of
meridians of the standard rotational hypersurface. They are called meridian
surfaces. We show that a meridian surface has a harmonic Gauss map if and
only if it is part of a plane. Further, we give necessary and sufficient
conditions for a meridian surface to have pointwise 1-type Gauss map and
find all meridian surfaces with pointwise 1-type Gauss map.
\end{abstract}

\maketitle

\section{Introduction}

The study of submanifolds of Euclidean space or pseudo-Euclidean space via
the notion of finite type immersions began in the late 1970's with the
papers \cite{Ch1,Ch2} of B.-Y. Chen and has been extensively carried out
since then. An isometric immersion $x:M$ $\rightarrow $ $\mathbb{E}^{m}$ of
a submanifold $M$ in Euclidean $m$-space $\mathbb{E}^{m}$ is said to be of
\emph{finite type} \cite{Ch1} if $x$ identified with the position vector
field of $M$ in $\mathbb{E}^{m}$ can be expressed as a finite sum of
eigenvectors of the Laplacian $\Delta $ of $M$, i.e.
\begin{equation*}
x=x_{0}+\sum_{i=1}^{k}x_{i},
\end{equation*}%
where $x_{0}$ is a constant map, $x_{1},x_{2},...,x_{k}$ are non-constant
maps such that $\Delta x_i=\lambda _{i}x_{i},$ $\lambda _{i}\in \mathbb{R}$,
$1\leq i\leq k.$ If $\lambda _{1},\lambda _{2},...,\lambda _{k}$ are
different, then $M$ is said to be of \emph{$k$-type}. Many results on finite
type immersions have been collected in the survey paper \cite{Ch3}.
Similarly, a smooth map $\phi $ of an $n$-dimensional Riemannian manifold $M$
of $\mathbb{E}^{m}$ is said to be of finite type if $\phi $ is a finite sum
of $\mathbb{E}^{m}$-valued eigenfunctions of $\Delta $. The notion of finite
type immersion is naturally extended to the Gauss map $G$ on $M$ in
Euclidean space \cite{CP}. Thus, a submanifold $M$ of Euclidean space has
\emph{1-type Gauss map} $G$, if $G$ satisfies $\Delta G=\mu (G+C)$ for some $%
\mu \in \mathbb{R}$ and some constant vector $C$ (of. \cite{BB}, \cite{BCV},
\cite{BV}, \cite{KY1}). However, the Laplacian of the Gauss map of some
typical well-known surfaces such as the helicoid, the catenoid and the right
cone in the Euclidean 3-space $\mathbb{E}^{3}$ takes a somewhat different
form, namely, $\Delta G=\lambda (G+C)$ for some non-constant function $%
\lambda $ and some constant vector $C$. Therefore, it is worth studying the
class of surfaces satisfying such an equation. A submanifold $M$ of the
Euclidean space $\mathbb{E}^{m}$ is said to have \emph{pointwise 1-type
Gauss map} if its Gauss map $G$ satisfies%
\begin{equation}
\Delta G=\lambda (G+C)  \label{A1}
\end{equation}%
for some non-zero smooth function $\lambda $ on $M$ and some constant vector
$C$ \cite{CK}. A pointwise 1-type Gauss map is called \emph{proper} if the
function $\lambda $ defined by (\ref{A1}) is non-constant. A submanifold
with pointwise 1-type Gauss map is said to be of the \emph{first kind} if
the vector $C$ in (\ref{A1}) is zero. Otherwise, the pointwise 1-type Gauss
map is said to be of the \emph{second kind} (\cite{CCK}, \cite{CK}, \cite%
{KY2}, \cite{KY3}). In \cite{CK} M. Choi and Y. Kim characterized the
minimal helicoid in terms of pointwise 1-type Gauss map of the first kind.
Also, together with B. Y. Chen, they proved that surfaces of revolution with
pointwise 1-type Gauss map of the first kind coincide with surfaces of
revolution with constant mean curvature \cite{CCK}. Moreover, they
characterized the rational surfaces of revolution with pointwise 1-type
Gauss map. In \cite{Yo2} D. Yoon studied Vranceanu rotation surfaces in
Euclidean 4-space $\mathbb{E}^{4}.$ He obtained classification theorems for
the flat Vranceanu rotation surfaces with 1-type Gauss map and an equation
in terms of the mean curvature vector \cite{Yo1}. For the general case see  \cite{A}.

The study of meridian surfaces in the Euclidean 4-space $\mathbb{E}^{4}$ was
first introduced by G. Ganchev and the third author in \cite{GM1}. The
meridian surfaces are one-parameter systems of meridians of the standard
rotational hypersurface in $\mathbb{E}^{4}$. In this paper we investigate
the meridian surfaces with pointwise 1-type Gauss map. We give necessary and
sufficient conditions for a meridian surface to have pointwise 1-type Gauss
map and find all meridian surfaces with pointwise 1-type Gauss map of first
and second kind.

\vskip 5mm

\section{ Preliminaries}

In the present section we recall definitions and results of \cite{Ch0}. Let $%
x:M$ $\rightarrow $ $\mathbb{E}^{m}$ be an immersion from an n-dimensional
connected Riemannian manifold $M$ into an $m-$dimensional Euclidean space $%
\mathbb{E}^{m}.$ We denote by $\langle , \rangle$ the metric tensor of $%
\mathbb{E}^{m}$ as well as the induced metric on M. Let $\nabla ^{\prime }$
be the Levi-Civita connection of $\mathbb{E}^{m}$ and $\nabla $ the induced
connection on $M$. Then the Gauss and Weingarten formulas are given
respectively by%
\begin{equation}
\nabla ^{\prime }_{X}Y=\nabla _{X}Y+h(X,Y),  \notag  \label{B1}
\end{equation}%
\begin{equation}
\nabla ^{\prime }_{X}\xi =-A_{\xi }X+D_{X}\xi  \notag  \label{B2}
\end{equation}%
where $X,Y$ are vector fields tangent to $M$ and $\xi $ is a vector field
normal to $M$. Moreover, $h$ is the second fundamental form, $D$ is the
linear connection induced in the normal bundle $T^{\perp }M$, called normal
connection, and $A_{\xi }$ is the shape operator in the direction of $\xi $
that is related with $h$ by
\begin{equation*}
\langle h(X,Y),\xi \rangle= \langle A_{\xi }X,Y \rangle.
\end{equation*}

The covariant differentiation $\overline{\nabla }h$ of the second
fundamental form $h$ on the direct sum of the tangent bundle and the normal
bundle $TM\oplus T^{\perp }M$ of $M$ is defined by
\begin{equation*}
(\overline{\nabla }_{X}h)(Y,Z)=D_{X}h(Y,Z)-h(\nabla _{X}Y,Z)-h(Y,\nabla
_{X}Z)
\end{equation*}%
for any vector fields $X,Y$ and $Z$ tangent to $M$. The Codazzi equation is
given by
\begin{equation}
(\overline{\nabla }_{X}h)(Y,Z)=(\overline{\nabla }_{Y}h)(X,Z).  \notag
\label{B3}
\end{equation}%
We denote by $R$ the curvature tensor associated with $\nabla$, i.e.
\begin{equation}
R(X,Y)Z=\nabla _{X}\nabla _{Y}Z-\nabla _{Y}\nabla _{X}Z-\nabla _{\lbrack
X,Y]}Z.  \notag  \label{B4}
\end{equation}
The equations of Gauss and Ricci are given, respectively, by
\begin{eqnarray}
& & \langle R(X,Y)Z,W \rangle = \langle h(X,W),h(Y,Z) \rangle - \langle
h(X,Z),h(Y,W) \rangle,  \notag  \label{B5} \\
& & \langle R^{\bot }(X,Y)\xi ,\eta \rangle = \langle [A_{\xi },A\eta ]X,Y
\rangle,  \notag  \label{B6}
\end{eqnarray}%
for vector fields $X,Y,Z,W$ tangent to $M$ and $\xi ,\eta $ normal to $M$.

The mean curvature vector field $H$ of an $n$-dimensional submanifold $M$ in
$\mathbb{E}^{m}$ is given by
\begin{equation*}
H = \frac{1}{n}\, trace\, h.
\end{equation*}

A submanifold $M$ is said to be minimal (respectively, totally geodesic) if $%
H \equiv 0$ (respectively, $h\equiv 0$).

We shall recall the definition of Gauss map $G$ of a submanifold $M$. Let $%
G(n,m)$ denote the Grassmannian manifold consisting of all oriented $n$%
-planes through the origin of $\mathbb{E}^{m}$ and $\wedge ^{n}\mathbb{E}%
^{m} $ be the vector space obtained by the exterior product of $n$ vectors
in $\mathbb{E}^{m}.$ In a natural way, we can identify $\wedge ^{n}\mathbb{E}%
^{m} $ with some Euclidean space $\mathbb{E}^{N}$ where $N=\left(
\begin{array}{c}
m \\
n%
\end{array}%
\right) .$ Let $\left \{ e_{1},\dots,e_{n},e_{n+1},\dots,e_{m}\right \} $ be
an adapted local orthonormal frame field in $\mathbb{E}^{m}$ such that $%
e_{1},e_{2},\dots,$ $e_{n},$ are tangent to $M$ and $e_{n+1}, e_{n+2},
\dots, e_{m}$ are normal to $M$. The map $G:M\rightarrow G(n,m)$ defined by $%
G(p)=(e_{1}\wedge e_{2}\wedge \dots \wedge $ $e_{n})(p)$ is called the Gauss
map of $M$. It is a a smooth map which carries a point $p$ in $M$ into the
oriented $n$-plane in $\mathbb{E}^{m}$ obtained by the parallel translation
of the tangent space of $M$ at $p$ in $\mathbb{E}^{m}.$

For any real function $\phi $ on $M$ the Laplacian of $\phi $ is defined by

\begin{equation}
\Delta \phi =-\sum_{i}(\nabla ^{\prime }_{e_{i}}\nabla ^{\prime
}_{e_{i}}\phi -\nabla ^{\prime }_{\nabla _{e_{i}}e_{i}}\phi ).  \label{B7}
\end{equation}

\section{Classification of meridian surfaces with pointwise 1-type Gauss map}

Let $\{e_{1},e_{2},e_{3},e_{4}\}$ be the standard orthonormal frame in $%
\mathbb{E}^{4}$, and $S^{2}(1)$ be the 2-dimensional sphere in $\mathbb{E}%
^{3}=span\{e_{1},e_{2},e_{3}\}$, centered at the origin $O$. We consider a
smooth curve $c:r=r(v),\,v\in J,\, \,J\subset \mathbb{R}$ on $S^{2}(1)$,
parameterized by the arc-length ($r^{\prime }{}^{2}(v)=1$). Let $%
t(v)=r^{\prime }(v)$ be the tangent vector field of $c$. We consider the
moving frame field $\{t(v),n(v),r(v)\}$ of the curve $c$ on $S^{2}(1)$. With
respect to this orthonormal frame field the following Frenet formulas hold:
\begin{equation}
\begin{array}{l}
\vspace{2mm}r^{\prime }=t; \\
\vspace{2mm}t^{\prime }=\kappa \,n-r; \\
\vspace{2mm}n^{\prime }=-\kappa \,t,%
\end{array}
\label{C1}
\end{equation}%
where $\kappa(v) = \langle t^{\prime }(v),n(v) \rangle$ is the spherical
curvature of $c$.

Let $f=f(u),\, \,g=g(u)$ be non-zero smooth functions, defined in an
interval $I\subset \mathbb{R}$, such that $(f^{\prime }(u))^{2}+(g^{\prime
}(u))^{2}=1,\, \,u\in I$. We consider the surface $M^{2}$ in $\mathbb{E}^{4}$
constructed in the following way:
\begin{equation}
M^{2}:z(u,v)=f(u)\,r(v)+g(u)\,e_{4},\quad u\in I,\,v\in J  \label{C2}
\end{equation}%
(see \cite{GM1}).

The surface $M^{2}$ lies on the rotational hypersurface $M^{3}$ in $\mathbb{E%
}^{4}$ obtained by the rotation of the meridian curve $\alpha :u\rightarrow
(f(u),g(u))$ about the $Oe_{4}$-axis in $\mathbb{E}^{4}$. $M^{2}$ is called
a \textit{meridian surface} on $M^3$ since it is a one-parameter system of
meridians of $M^{3}$.

The tangent space of $M^{2}$ is spanned by the vector fields:
\begin{equation}  \label{C3}
\begin{array}{l}
\vspace{2mm}z_{u}=f^{\prime }r+g^{\prime }e_{4}; \\
\vspace{2mm}z_{v}=f\,t,%
\end{array}%
\end{equation}%
and hence, the coefficients of the first fundamental form of $M^{2}$ are $%
E=1;\, \,F=0;\, \,G=f^{2}(u)$. Taking into account \eqref{C1} and \eqref{C3}%
, we calculate the second partial derivatives of $z(u,v)$: \
\begin{equation*}
\begin{array}{l}
\vspace{2mm}z_{uu}=f^{\prime \prime }r+g^{\prime \prime }\,e_{4}; \\
\vspace{2mm}z_{uv}=f^{\prime }t; \\
\vspace{2mm}z_{vv}=f\kappa \,n-f\,r.%
\end{array}%
\end{equation*}%
Let us denote $x=z_{u},\, \,y=\displaystyle{\frac{z_v}{f}=t}$ and consider
the following orthonormal normal frame field of $M^{2}$:
\begin{equation}
n_{1}=n(v); \quad n_{2}=-g^{\prime }(u)\,r(v)+f^{\prime }(u)\,e_{4}.  \notag
\label{C4}
\end{equation}%
Thus we obtain a positive orthonormal frame field $\{x,y,n_{1},n_{2}\}$ of $%
M^{2}$. We denote by $\kappa _{\alpha }$ the curvature of the meridian curve
$\alpha $, i.e.
\begin{equation}
\kappa _{\alpha }(u)=f^{\prime }(u)\,g^{\prime \prime }(u)-g^{\prime
}(u)f^{\prime \prime }(u).  \notag  \label{C5}
\end{equation}
By covariant differentiation with respect to $x$ and $y$, and a
straightforward calculation we obtain
\begin{equation}
\begin{array}{l}
\vspace{2mm} \nabla ^{\prime }_{x}x=\kappa _{\alpha }\,n_{2}; \\
\vspace{2mm} \nabla ^{\prime }_{x}y=0; \\
\vspace{2mm} \nabla ^{\prime }_{y}x=\displaystyle{\frac{f^{\prime }}{f}\,y};
\\
\vspace{2mm} \nabla ^{\prime }_{y}y=\displaystyle{-\frac{f^{\prime }}{f}\,x+%
\frac{\kappa }{f}\,n_{1}+\frac{g^{\prime }}{f}\,n_{2}};%
\end{array}
\label{C6}
\end{equation}
and
\begin{equation}  \label{C7}
\begin{array}{l}
\vspace{2mm} \nabla ^{\prime }_{x}n_{1}=0; \\
\vspace{2mm} \nabla ^{\prime }_{y}n_{1}=\displaystyle{-\frac{\kappa }{f}\,y};
\\
\vspace{2mm} \nabla ^{\prime }_{x}n_{2}=-\kappa _{\alpha }\,x; \\
\vspace{2mm} \nabla ^{\prime }_{y}n_{2}=\displaystyle{-\frac{g^{\prime }}{f}%
\,y}.%
\end{array}%
\end{equation}
where $\kappa(v)$ and $\kappa _{\alpha}(u)$ are the curvatures of the
spherical $c$ and the meridian curve $\alpha$, respectively (see \cite{GM1}).

Equalities (\ref{C7}) imply the following result.

\begin{lem}
Let $M^{2}$ be a meridian surface given with the surface patch \eqref{C2}.
Then
\begin{equation}
A_{n_1}=\left[
\begin{array}{ll}
0 & 0 \\
0 & \displaystyle{\frac{\kappa}{f}}%
\end{array}
\right], \qquad A_{n_2}=\left[
\begin{array}{ll}
\kappa _{\alpha } & 0 \\
0 & \displaystyle{\frac{g^{\prime }}{f}}%
\end{array}
\right] .  \notag  \label{C8}
\end{equation}
\end{lem}

\vskip 3mm So, the Gauss curvature is given by
\begin{eqnarray}
K &= &{\frac{\kappa _{\alpha} \,g^{\prime }}{f}}  \notag  \label{C9}
\end{eqnarray}
and the mean curvature vector field $H$ of $M^{2}$ is
\begin{eqnarray}  \label{C10}
H &=&{\frac{\kappa }{2f} \,n_1+\frac{\kappa _{\alpha }f+g^{\prime }}{2f}
\,n_2.}  \notag
\end{eqnarray}

The Gauss map $G$ of $M^{2}$ is defined by $G=x\wedge y$. Using (\ref{B7}), (%
\ref{C6}), and (\ref{C7}) we calculate that the Laplacian of the Gauss map
is expressed as
\begin{eqnarray}
\Delta G &=& {\frac{(f \kappa_{\alpha})^{2}+\kappa ^{2} + g^{\prime 2}}{f^{2}%
}}\, x\wedge y - \frac{\kappa^{\prime }}{f^{2}}\, x\wedge n_{1}  \label{C11}
\\
&&- {\frac{\kappa f^{\prime }}{f^{2}}}\, y\wedge n_{1} - \frac{f^{\prime}
g^{\prime} - f ( f \kappa_{\alpha} )^{\prime} }{f^{2}}\, y\wedge n_{2},
\notag
\end{eqnarray}%
where $\kappa^{\prime }= \displaystyle{\frac{d}{dv}(\kappa)}$.

First, we suppose that the Gauss map of $M^2$ is harmonic, i.e. $\Delta G =
0 $. Then from \eqref{C11} we get
\begin{eqnarray}
\kappa _{\alpha } &=&0;  \notag \\
\kappa &=&0;  \label{C7*} \\
g^{\prime } &=&0.  \notag
\end{eqnarray}

So, \eqref{C6} and \eqref{C7*} imply that $M^2$ is a totally geodesic
surface in $\mathbb{E}^{4}$. Conversely, if $M^2$ is totally geodesic, then $%
\Delta G = 0$.

Thus we obtain the following result.

\begin{thm}
Let $M^2$ be a meridian surfaces in the Euclidean space $\mathbb{E}^{4}$.
The Gauss map of $M^2$ is harmonic if and only if $M^2$ is part of a plane.
\end{thm}

Now, we suppose that the meridian surface $M^2$ is of pointwise 1-type Gauss
map, i.e. $G$ satisfies \eqref{A1}, where $\lambda \neq 0$. Then, from
equalities \eqref{A1} and \eqref{C11} we get
\begin{eqnarray}
\lambda +\lambda \; \langle C,x\wedge y \rangle & = & \displaystyle{{\frac{%
(f \kappa_{\alpha})^{2}+\kappa ^{2}+g^{\prime 2}}{f^{2}}}};  \notag \\
\lambda \; \langle C,x\wedge n_{1}\rangle &=& -\displaystyle{\frac{%
\kappa^{\prime }}{f^{2}}};  \label{C12} \\
\lambda \; \langle C,y\wedge n_{1}\rangle &=& -\displaystyle{{\frac{\kappa
f^{\prime }}{f^{2}}}};  \notag \\
\lambda \; \langle C,y\wedge n_{2}\rangle &=& -\displaystyle{\frac{%
f^{\prime} g^{\prime} - f ( f \kappa_{\alpha} )^{\prime} }{f^{2}}}.  \notag
\end{eqnarray}
Using \eqref{C11} we obtain
\begin{eqnarray}
\lambda \; \langle C,x\wedge n_{2}\rangle=0;  \notag \\
\lambda \; \langle C,n_{1}\wedge n_{2}\rangle=0.  \label{C13}
\end{eqnarray}

Differentiating (\ref{C13}) with respect to $u$ and $v$ we get%
\begin{eqnarray}
\kappa _{\alpha } \; \langle C,x\wedge n_{1}\rangle &= &0;  \notag \\
\frac{f{}^{\prime }}{f{}}\; \langle C,y\wedge n_{2}\rangle-\frac{g{}^{\prime
}}{f{}} \; \langle C,x\wedge y\rangle & =&0;  \label{C16} \\
-\frac{\kappa }{f{}}\; \langle C,y\wedge n_{2}\rangle +\frac{g{}^{\prime }}{%
f{}} \; \langle C,y\wedge n_{1}\rangle &=&0.  \notag
\end{eqnarray}
Since $\lambda \neq 0$ equalities \eqref{C12} and \eqref{C16} imply
\begin{equation}  \label{C18a}
\begin{array}{l}
\vspace{2mm} \kappa _{\alpha } \kappa^{\prime }=0; \\
\vspace{2mm} \kappa \, (f \kappa_{\alpha})^{\prime} = 0; \\
\vspace{2mm} \lambda f^2 g^{\prime} = g^{\prime }\left(1 + (f
\kappa_{\alpha})^{2}+\kappa^{2}\right) - f f^{\prime }( f \kappa_{\alpha}
)^{\prime }.%
\end{array}%
\end{equation}

We distinguish the following cases.

\vskip 2mm \textbf{Case I:} $g^{\prime }= 0$. In such case $\kappa_{\alpha}
= 0$. Then equality \eqref{C11} implies that
\begin{equation}  \label{C11a}
\Delta G = {\frac{\kappa ^{2}}{f^{2}}}\, x\wedge y - \frac{\kappa^{\prime }}{%
f^{2}}\, x\wedge n_{1} - {\frac{\kappa f^{\prime }}{f^{2}}}\, y\wedge n_{1}.
\end{equation}

If we assume that $M^2$ has pointwise 1-type Gauss map of the first kind,
i.e. $C = 0$, then from \eqref{C11a} we get $\kappa^{\prime }= 0$ and $%
\kappa f^{\prime }= 0$, which imply $\kappa = 0$ since $f^{\prime }\neq 0$.
Hence $\Delta G = 0$, which contradicts the assumption that $\lambda \neq 0$%
. Consequently, in the case $g^{\prime }= 0$ there are no meridian surfaces
of pointwise 1-type Gauss map of the first kind.

Now we consider meridian surfaces of pointwise 1-type Gauss map of the
second kind, i.e. $C \neq 0$. So we suppose that $\kappa \neq 0$. From
equalities \eqref{A1} and \eqref{C11a} we obtain
\begin{equation}  \label{C30}
C = \left( {\frac{\kappa ^{2}}{\lambda f^{2}}} -1 \right)\, x\wedge y -
\frac{\kappa^{\prime }}{\lambda f^{2}}\, x\wedge n_{1} - {\frac{\kappa
f^{\prime }}{\lambda f^{2}}}\, y\wedge n_{1}.
\end{equation}
Using \eqref{C6}, \eqref{C7} and \eqref{C30} we obtain
\begin{equation*}
\begin{array}{lll}
\vspace{2mm} \nabla^{\prime }_xC & = & \displaystyle{\kappa ^{2} \left( {%
\frac{1}{\lambda f^{2}}} \right)^{\prime }_u \, x\wedge y - \kappa^{\prime
}\left( {\frac{1}{\lambda f^{2}}} \right)^{\prime }_u\, x\wedge n_{1} -
\kappa f^{\prime} \left( {\frac{1}{\lambda f^{2}}} \right)^{\prime }_u \,
y\wedge n_{1};} \\
\vspace{2mm} \nabla^{\prime }_yC & = & \displaystyle{\frac{\kappa}{\lambda^2
f^3} \left( 3 \kappa^{\prime }\lambda - \kappa \lambda^{\prime }_v \right)
\, x\wedge y + \frac{1}{\lambda^2 f^3} \left(- \kappa^{\prime \prime
}\lambda + k^{\prime }\lambda^{\prime }_v + \kappa^3 \lambda + \kappa
\lambda - \kappa \lambda^2 f^2 \right)\, x\wedge n_{1}} \\
&  & \displaystyle{+ \frac{f^{\prime }}{\lambda^2 f^3} \left(- 2
\kappa^{\prime }\lambda + \kappa \lambda^{\prime }_v \right) \, y\wedge n_{1}%
}.%
\end{array}%
\end{equation*}
The last formulas imply that $C = const$ if and only if $\kappa = const$ and
$\lambda = \displaystyle{\frac{\kappa^2 + 1}{f^2}}$.

The condition $\kappa = const \neq 0$ implies that the curve $c$ on $%
S^{2}(1) $ is a circle with non-zero constant spherical curvature. Since $%
g^{\prime }= 0$ and $(f^{\prime 2 }+ (g^{\prime 2 }= 1$ we get $f(u) = \pm u
+ a$, $g(u) = b$, where $a=const$, $b=const$. In this case $M^2$ is a
developable ruled surface. Moreover, from \eqref{C7} it follows that $%
\nabla^{\prime }_x n_2 = 0; \, \, \nabla^{\prime }_y n_2 =0$, which implies
that $M^2$ lies in the 3-dimensional space spanned by $\{x, y, n_1\}$.

Conversely, if $g^{\prime }= 0$ and $\kappa = const$, by direct computation
we get
\begin{equation*}
\Delta G = {\frac{\kappa ^{2} + 1}{f^{2}}}(G + C),
\end{equation*}
where $C = \displaystyle{{-\frac{1}{\kappa ^{2} + 1}}\, x\wedge y - \frac{%
\kappa f^{\prime }}{\kappa ^{2} + 1}\, y\wedge n_{1}}$. Hence, $M^2$ is a
surface with pointwise 1-type Gauss map of the second kind.

Summing up we obtain the following result.

\begin{thm}
\label{T:second kind-0} Let $M^2$ be a meridian surface given with
parametrization \eqref{C2} and $g^{\prime }=0$. Then $M^2$ has pointwise
1-type Gauss map of the second kind if and only if the curve $c$ is a circle
with non-zero constant spherical curvature and the meridian curve $\alpha$
is determined by $f(u) = \pm u + a; \, \, g(u) = b$, where $a=const$, $%
b=const $. In this case $M^2$ is a developable ruled surface lying in
3-dimensional space.
\end{thm}

\vskip 2mm \textbf{Case II:} $g^{\prime }\neq 0$. In such case from the
third equality of \eqref{C18a} we obtain
\begin{equation}  \label{C34}
\lambda =\frac{g^{\prime }\left(1 + (f
\kappa_{\alpha})^{2}+\kappa^{2}\right) - f f^{\prime }( f \kappa_{\alpha}
)^{\prime }}{f^2 g^{\prime }}.
\end{equation}

First we shall consider the case of pointwise 1-type Gauss map surfaces of
the first kind. From \eqref{C11} it follows that $M^2$ is of the first kind (%
$C = 0$) if and only if
\begin{equation}  \label{C31}
\begin{array}{l}
\vspace{2mm} \kappa^{\prime }= 0; \\
\vspace{2mm} \kappa f^{\prime }= 0; \\
\vspace{2mm} f^{\prime }g^{\prime }- f ( f \kappa_{\alpha} )^{\prime }= 0.%
\end{array}%
\end{equation}

The first equality of \eqref{C31} implies that $\kappa = const$. There are
two subcases:

\vskip 2mm \hskip 10mm 1. $\kappa = 0$. Then the meridian curve $\alpha$ is
determined by the equation
\begin{equation}  \label{C32}
f^{\prime }g^{\prime }- f ( f \kappa_{\alpha} )^{\prime }= 0.
\end{equation}
The equalities $\kappa_{\alpha} = f^{\prime }g^{\prime \prime }- g^{\prime
}f^{\prime \prime }$ and $f^{\prime 2 }+ g^{\prime 2 }= 1$ imply that $%
\kappa _{\alpha } = \displaystyle{-{\frac{f^{\prime \prime }}{g^{\prime }}}}$%
. Hence equation \eqref{C32} can be rewritten in the form
\begin{equation}  \label{C33}
f^{\prime }\sqrt{1 - f^{\prime 2}} + f \left( \frac{f f^{\prime \prime }}{%
\sqrt{1 - f^{\prime 2}}}\right)^{\prime }= 0.
\end{equation}
Since $\kappa = 0$, $M^2$ lies in the 3-dimensional space spanned by $\{x,
y, n_2\}$.

Conversely, if $\kappa = 0$ and the meridian curve $\alpha$ is determined by
a solution $f(u)$ of differential equation \eqref{C33}, the function $g(u)$
is defined by $g^{\prime }=\sqrt{1-f^{\prime 2}}$, then the surface $M^2$,
parameterized by \eqref{C2}, is a surface of pointwise 1-type Gauss map of
the first kind.

\vskip 2mm \hskip 10mm 2. $\kappa \neq 0$. Then the second equality of %
\eqref{C31} implies that $f^{\prime }= 0$. In this case $f(u) = a; \, \,
g(u) = \pm u + b$, where $a = const$, $b = const$. By a result of \cite{GM1}%
, $M^2 $ is a developable ruled surface in a 3-dimensional space, since $%
\kappa_{\alpha} = 0$ and $\kappa = const$. It follows from \eqref{C34} that $%
\lambda = \displaystyle{\frac{1 + \kappa^2}{a^2}} = const$, which implies
that $M^2$ has 1-type Gauss map, i.e. $M^2$ is non-proper. The converse is
also true.

\vskip 2mm Thus we obtain the following result.

\begin{thm}
\label{T:first kind} Let $M^2$ be a meridian surface given with
parametrization \eqref{C2} and $g^{\prime }\neq 0$. Then $M^2$ has pointwise
1-type Gauss map of the first kind if and only if one of the following holds:

(i) the curve $c$ is a great circle on $S^{2}(1)$ and the meridian curve $%
\alpha$ is determined by the solutions of the following differential
equation
\begin{equation*}
f^{\prime }\sqrt{1 - f^{\prime 2}} + f \left( \frac{f f^{\prime \prime }}{%
\sqrt{1 - f^{\prime 2}}}\right)^{\prime }= 0;
\end{equation*}

(ii) the curve $c$ is a circle on $S^2(1)$ with non-zero constant spherical
curvature and the meridian curve $\alpha$ is determined by $f(u) = a; \, \,
g(u) = \pm u + b$, where $a=const$, $b=const$. In this case $M^2$ is a
developable ruled surface in a 3-dimensional space. Moreover, $M^2$ is
non-proper.
\end{thm}

\vskip 2mm Now we shall consider the case of pointwise 1-type Gauss map
surfaces of the second kind. It follows from equalities \eqref{C18a} that
there are three subcases.

\vskip 2mm \hskip 10mm 1. $\kappa_{\alpha} = 0$. In this subcase
\begin{equation}  \label{C11b}
\Delta G = {\frac{\kappa^2 + g^{\prime 2}}{f^{2}}}\, x\wedge y - \frac{%
\kappa^{\prime }}{f^{2}}\, x\wedge n_{1} - \frac{\kappa f^{\prime }}{f^2}\,
y\wedge n_{1} - \frac{f^{\prime }g^{\prime }}{f^2}\, y\wedge n_{2}.
\end{equation}
From equalities \eqref{A1} and \eqref{C11b} we obtain
\begin{equation*}
C = \left( {\frac{\kappa ^{2} + g^{\prime 2}}{\lambda f^{2}}} -1 \right)\,
x\wedge y - \frac{\kappa^{\prime }}{\lambda f^{2}}\, x\wedge n_{1} - \frac{%
\kappa f^{\prime }}{\lambda f^2}\, y\wedge n_{1} - \frac{f^{\prime
}g^{\prime }}{\lambda f^2}\, y\wedge n_{2}.
\end{equation*}
The third equality in \eqref{C18a} implies that in this case $\lambda = %
\displaystyle{\frac{1 + \kappa^2}{f^2}}$ and hence, $C$ is expressed as
follows:
\begin{equation}  \label{C36}
C = - \frac{1}{1 + \kappa^2} \left( f^{\prime 2}\, x\wedge y +
\kappa^{\prime }\, x\wedge n_{1} + \kappa f^{\prime }\, y\wedge n_{1} +
f^{\prime }g^{\prime }\, y\wedge n_{2} \right).
\end{equation}

Using \eqref{C6}, \eqref{C7} and \eqref{C36} we obtain
\begin{equation*}
\begin{array}{lll}
\vspace{2mm}\nabla _{x}^{\prime }C & = & -\displaystyle{\frac{1}{1+\kappa
^{2}}}\left( 2f^{\prime }f^{\prime \prime }\,x\wedge y+\kappa f^{\prime
\prime }\,y\wedge n_{1}+(f^{\prime }g^{\prime \prime }+f^{\prime \prime
}g^{\prime })\,y\wedge n_{2}\right) ; \\
\vspace{2mm}\nabla _{y}^{\prime }C & = & \displaystyle{\frac{1}{f(1+\kappa
^{2})^{2}}\left( \left( 2\kappa \kappa ^{\prime }f^{\prime }{}^{2}+\kappa
\kappa ^{\prime }(1+\kappa ^{2})\right) \,x\wedge y+\left( 2\kappa \kappa
^{\prime }{}^{2}-(1+\kappa ^{2})\kappa ^{\prime \prime }\right) \,x\wedge
n_{1}\right) } \\
&  & \displaystyle{\ +\frac{1}{f(1+\kappa ^{2})^{2}}\left(-2\kappa^{\prime }f^{\prime }\,y\wedge n_{1}+2\kappa \kappa ^{\prime
}f^{\prime }g^{\prime
}\,y\wedge n_{2}\right) }.%
\end{array}%
\end{equation*}

The last formulas imply that $C = const$ if and only if $\kappa = const$, $%
f^{\prime }= a = const$, $g^{\prime }= b = const$, $a^2 + b^2 = 1$.

The condition $\kappa = const$ implies that the curve $c$ is a circle on $%
S^{2}(1)$. The meridian curve $\alpha$ is given by $f(u) = a u + a_1; \, \,
g(u) = b u + b_1$, where $a_1=const$, $b_1=const$. In this case $M^2$ is a
developable ruled surface lying in a 3-dimensional space.

Conversely, if $f(u) = a u + a_1; \, \,g(u) = b u + b_1$ and $\kappa = const$%
, then
\begin{equation*}
\Delta G = {\frac{\kappa^2 + b^2}{f^{2}}}\, x\wedge y - \frac{\kappa a}{f^2}%
\, y\wedge n_{1} - \frac{ab }{f^2}\, y\wedge n_{2}.
\end{equation*}
Hence, by direct computation we get
\begin{equation*}
\Delta G = {\frac{1 + \kappa ^{2}}{f^{2}}}(G + C),
\end{equation*}
where $C = \displaystyle{- \frac{a}{1 + \kappa^2} \left( a\, x\wedge y +
\kappa \, y\wedge n_{1} + b \, y\wedge n_{2} \right)}$. Consequently, $M^2$
is a surface of pointwise 1-type Gauss map of the second kind.

\vskip 2mm \hskip 10mm 2. $\kappa = 0$. In this subcase
\begin{equation}  \label{C11c}
\Delta G = {\frac{(f \kappa_{\alpha})^{2} +g^{\prime 2}}{f^{2}}}\, x\wedge y
- \frac{f^{\prime} g^{\prime} - f ( f \kappa_{\alpha} )^{\prime} }{f^2}\,
y\wedge n_{2}.
\end{equation}

From equalities \eqref{A1} and \eqref{C11c} we obtain
\begin{equation*}
C = \left( \frac{(f \kappa_{\alpha})^{2} +g^{\prime 2}}{\lambda f^{2}} -1
\right)\, x\wedge y - \frac{f^{\prime} g^{\prime} - f ( f \kappa_{\alpha}
)^{\prime} }{\lambda f^2}\, y\wedge n_{2}.
\end{equation*}
Using the third equality of \eqref{C18a} we obtain that $C$ is expressed as
follows:
\begin{equation}  \label{C37}
C = - \frac{f^{\prime} g^{\prime} - f ( f \kappa_{\alpha} )^{\prime}}{%
\lambda f^2} \left( \frac{f^{\prime }}{g^{\prime }} \, x\wedge y + y\wedge
n_{2} \right),
\end{equation}
where $\lambda = \displaystyle{\frac{1}{f^2} \left(1 + (f
\kappa_{\alpha})^{2} - \frac{f f^{\prime }}{g^{\prime }}( f \kappa_{\alpha}
)^{\prime}\right)}$. We denote
\begin{equation}  \label{C37-a}
\varphi = \displaystyle{- \frac{f^{\prime} g^{\prime} - f ( f
\kappa_{\alpha} )^{\prime}}{\lambda f^2}}.
\end{equation}
Then equalities \eqref{C6}, \eqref{C7} and \eqref{C37} imply

\begin{equation}
\begin{array}{l}
\vspace{2mm}\nabla _{x}^{\prime }C=\displaystyle{\left( \left( \varphi \,%
\frac{f^{\prime }}{g^{\prime }}\right) ^{\prime }+\varphi \kappa _{\alpha
}\right) \,x\wedge y+\left( \varphi ^{\prime }-\varphi \, \frac{f^{\prime }}{%
g^{\prime }}\kappa _{\alpha }\right) \,y\wedge n_{2}}; \\
\vspace{2mm}\nabla _{y}^{\prime }C=0.%
\end{array}
\label{C37-b}
\end{equation}%
It follows from \eqref{C37-b} imply that $C=const$ if and only if $\varphi
^{\prime }=\displaystyle{\varphi \, \frac{f^{\prime }}{g^{\prime }}\kappa
_{\alpha }}$, or equivalently
\begin{equation}
\left( \ln \varphi \right) ^{\prime }=\displaystyle{\frac{f^{\prime }}{%
g^{\prime }}\kappa _{\alpha }}.  \label{C37-c}
\end{equation}%
Using that $f\kappa _{\alpha }=\displaystyle{-{\frac{ff^{\prime \prime }}{%
\sqrt{1-f^{\prime 2}}}}}$, from \eqref{C37-a} we get
\begin{equation}
\varphi =\frac{-\sqrt{1-f^{\prime 2}}\left( f(1-f^{\prime
2})(ff^{\prime \prime })' + f^2f^{\prime }f^{\prime \prime
2}+f^{\prime }(1-f^{\prime 2})^{2}\right) }{ff^{\prime
}(ff^{\prime \prime })^{\prime }(1-f^{\prime 2})+f^{2}f^{\prime
\prime 2}+(1-f^{\prime 2})^{2}}.  \label{C37-d}
\end{equation}%
Now, formulas \eqref{C37-c} and \eqref{C37-d} imply that $C=const$ if and
only if the function $f(u)$ is a solution of the following differential
equation
\begin{equation}
\left( \ln \frac{-\sqrt{1-f^{\prime 2}}\left( f(1-f^{\prime
2})(ff^{\prime \prime })' + f^2f^{\prime }f^{\prime \prime
}{}^{2}+f^{\prime }(1-f^{\prime 2})^{2}\right) }{ff^{\prime
}(ff^{\prime \prime })^{\prime }(1-f^{\prime 2})+f^{2}f^{\prime
\prime 2}+(1-f^{\prime 2})^{2}}\right) ^{\prime
}=-{\frac{f^{\prime }f^{\prime \prime }}{1-f^{\prime 2}}}.
\label{C37-eq}
\end{equation}

Conversely, if $\kappa = 0$ and the meridian curve $\alpha$ is determined by
a solution $f(u)$ of differential equation \eqref{C37-eq}, $g(u)$ is defined
by $g^{\prime }=\sqrt{1-f^{\prime 2}}$, then the surface $M^2$,
parameterized by \eqref{C2}, is a surface of pointwise 1-type Gauss map of
the second kind.

\vskip2mm \hskip10mm 3. $\kappa =const\neq 0$ and $f\kappa _{\alpha
}=a=const,\,a\neq 0$. In this subcase
\begin{equation}
\Delta G={\frac{a^{2}+\kappa ^{2}+g^{\prime 2}}{f^{2}}}\,x\wedge y-{\frac{%
\kappa f^{\prime }}{f^{2}}}\,y\wedge n_{1}-\frac{f^{\prime }g^{\prime }}{%
f^{2}}\,y\wedge n_{2}.  \label{C11d}
\end{equation}%
From equalities \eqref{A1}, \eqref{C34} and \eqref{C11d} we obtain
\begin{equation}
C=-\frac{1}{1+a^{2}+\kappa ^{2}}\left( f^{\prime 2}\,x\wedge y+\kappa
f^{\prime }\,y\wedge n_{1}+f^{\prime }g^{\prime }\,y\wedge n_{2}\right) .
\label{C38}
\end{equation}%
Then equalities \eqref{C6}, \eqref{C7} and \eqref{C38} imply

\begin{equation}  \label{C39}
\begin{array}{l}
\vspace{2mm} \nabla^{\prime }_x C = \displaystyle{- \frac{1}{1 + a^2 +
\kappa^2} \left( f^{\prime }f^{\prime \prime }\, x\wedge y + \kappa
f^{\prime \prime }\, y\wedge n_{1} + g^{\prime }f^{\prime \prime }\, y\wedge
n_{2} \right)}; \\
\vspace{2mm} \nabla^{\prime }_y C = 0.%
\end{array}%
\end{equation}

Formulas \eqref{C39} imply that $C = const$ if and only if $f^{\prime \prime
}= 0$. But, if $f^{\prime \prime }= 0$ then $\kappa_{\alpha} = 0$, which
contradicts the assumption that $f \kappa_{\alpha} \neq 0$.

Consequently, if $\kappa = const \neq 0$ and $f \kappa_{\alpha} = a = const,
\, a \neq 0$, then there are no meridian surfaces of pointwise 1-type Gauss
map of the second kind.

\vskip 2mm Summing up we obtain the following result.

\begin{thm}
\label{T:second kind} Let $M^2$ be a meridian surface given with
parametrization \eqref{C2} and $g^{\prime }\neq 0$. Then $M^2$ has pointwise
1-type Gauss map of the second kind if and only if one of the following
holds:

(i) the curve $c$ is a circle on $S^{2}(1)$ and the meridian curve $\alpha$
is determined by $f(u) = a u + a_1; \, \, g(u) = b u + b_1$, where $a$, $a_1$%
, $b$, $b_1$ are constants. In this case $M^2$ is a developable ruled
surface lying in a 3-dimensional space;

(ii) the curve $c$ is a great circle on $S^{2}(1)$ and the meridian curve $%
\alpha$ is determined by the solutions of the following differential
equation
\begin{equation*}
\left(\ln \frac{-\sqrt{1 - f^{\prime 2}} \left( f (1 - f^{\prime
2}) (f f^{\prime \prime })' + f^2f^{\prime }f^{\prime \prime 2 }+
f^{\prime }(1 - f^{\prime 2})^2 \right)}{f f^{\prime }(f f^{\prime
\prime })^{\prime
}(1 - f^{\prime 2}) + f^2 f^{\prime \prime 2 }+ (1 - f^{\prime 2})^2}%
\right)^{\prime }= -{\frac{f^{\prime }f^{\prime \prime }}{1 - f^{\prime 2}}}.
\end{equation*}
\end{thm}

\vskip 3mm Theorem \ref{T:second kind-0}, Theorem \ref{T:first kind}, and
Theorem \ref{T:second kind} describe all meridian surfaces with pointwise
1-type Gauss map.

\vskip 3mm \textbf{Acknowledgements:} This paper is prepared during the
third named author's visit to the Uluda\u{g} University, Bursa, Turkey in
January 2011.

\end{document}